\documentclass[12pt]{article}

\newcommand{\ol}{\setlength{\itemsep}{0pt.}\begin{enumerate}}
\newcommand{\eol}{\end{enumerate}\setlength{\itemsep}{-\parsep}}
\newcommand{\ignore}[1]{}
\def\R{{\bf R}}

\title{\bf Unleashing the power of Schrijver's permanental inequality with the help of the Bethe Approximation.
 }
 
\author{Leonid Gurvits \thanks{%
{\tt gurvits@lanl.gov}. Los Alamos National Laboratory, 
Los Alamos, NM. } 
}

\begin{document}



\maketitle

\begin{abstract}
Let $A \in \Omega_n$ be doubly-stochastic $n \times n$ matrix. Alexander Schrijver proved in 1998 the following remarkable inequality
\begin{equation} \label{le}
per(\widetilde{A}) \geq \prod_{1 \leq i,j \leq n} (1- A(i,j)); \widetilde{A}(i,j) =: A(i,j)(1-A(i,j)), 1 \leq i,j \leq n
\end{equation}
We prove in this paper the following generalization (or just clever reformulation) of (\ref{le}):\\
For all pairs of $n \times n$ matrices $(P,Q)$, where $P$ is nonnegative and $Q$ is doubly-stochastic
\begin{equation} \label{st}
\log(per(P)) \geq \sum_{1 \leq i,j \leq n} \log(1- Q(i,j)) (1- Q(i,j)) - \sum_{1 \leq i,j \leq n} Q(i,j) \log \left(\frac{Q(i,j)}{P(i,j)} \right
)
\end{equation}
The main co
rollary of (\ref{st}) is the following inequality for doubly-stochastic matrices:
$$
\frac{per(A)}{F(A)} \geq 1; F(A) =: \prod_{1 \leq i,j \leq n} \left(1- A(i,j)\right)^{1- A(i,j)}.
$$
{\bf We use this inequality to prove Friedland's conjecture on monomer-dimer entropy, so called {\it Asymptotic Lower Matching Conjecture}}\\
We present explicit doubly-stochastic $n \times n$ matrices $A$ with
the ratio $\frac{per(A)}{F(A)} = \sqrt{2}^{n}$ and conjecture that
$$
\max_{A \in \Omega_n}\frac{per(A)}{F(A)} \approx \left(\sqrt{2} \right)^{n}.
$$
If true, it would imply a deterministic poly-time algorithm to approximate the permanent of $n \times n$ nonnegative
matrices within the relative factor $\left(\sqrt{2} \right)^{n}$.\\

\end{abstract} 

 
\newtheorem{THEOREM}{Theorem}[section]
\newenvironment{theorem}{\begin{THEOREM} \hspace{-.85em} {\bf :} 
}%
                        {\end{THEOREM}}
\newtheorem{LEMMA}[THEOREM]{Lemma}
\newenvironment{lemma}{\begin{LEMMA} \hspace{-.85em} {\bf :} }%
                      {\end{LEMMA}}
\newtheorem{COROLLARY}[THEOREM]{Corollary}
\newenvironment{corollary}{\begin{COROLLARY} \hspace{-.85em} {\bf 
:} }%
                          {\end{COROLLARY}}
\newtheorem{PROPOSITION}[THEOREM]{Proposition}
\newenvironment{proposition}{\begin{PROPOSITION} \hspace{-.85em} 
{\bf :} }%
                            {\end{PROPOSITION}}
\newtheorem{DEFINITION}[THEOREM]{Definition}
\newenvironment{definition}{\begin{DEFINITION} \hspace{-.85em} {\bf 
:} \rm}%
                            {\end{DEFINITION}}
\newtheorem{EXAMPLE}[THEOREM]{Example}
\newenvironment{example}{\begin{EXAMPLE} \hspace{-.85em} {\bf :} 
\rm}%
                            {\end{EXAMPLE}}
\newtheorem{CONJECTURE}[THEOREM]{Conjecture}
\newenvironment{conjecture}{\begin{CONJECTURE} \hspace{-.85em} 
{\bf :} \rm}%
                            {\end{CONJECTURE}}
\newtheorem{PROBLEM}[THEOREM]{Problem}
\newenvironment{problem}{\begin{PROBLEM} \hspace{-.85em} {\bf :} 
\rm}%
                            {\end{PROBLEM}}
\newtheorem{QUESTION}[THEOREM]{Question}
\newenvironment{question}{\begin{QUESTION} \hspace{-.85em} {\bf :} 
\rm}%
                            {\end{QUESTION}}
\newtheorem{REMARK}[THEOREM]{Remark}
\newenvironment{remark}{\begin{REMARK} \hspace{-.85em} {\bf :} 
\rm}%
                            {\end{REMARK}}
\newtheorem{FACT}[THEOREM]{Fact}
\newenvironment{fact}{\begin{FACT} \hspace{-.85em} {\bf :} 
\rm}%
		            {\end{FACT}}

 
\newcommand{\thm}{\begin{theorem}}
\newcommand{\lem}{\begin{lemma}}
\newcommand{\pro}{\begin{proposition}}
\newcommand{\dfn}{\begin{definition}}
\newcommand{\rem}{\begin{remark}}
\newcommand{\xam}{\begin{example}}
\newcommand{\cnj}{\begin{conjecture}}
\newcommand{\prb}{\begin{problem}}
\newcommand{\que}{\begin{question}}
\newcommand{\cor}{\begin{corollary}}
\newcommand{\fac}{\begin{fact}}

\newcommand{\prf}{\noindent{\bf Proof:} }
\newcommand{\ethm}{\end{theorem}}
\newcommand{\elem}{\end{lemma}}
\newcommand{\epro}{\end{proposition}}
\newcommand{\edfn}{\bbox\end{definition}}
\newcommand{\erem}{\bbox\end{remark}}
\newcommand{\exam}{\bbox\end{example}}
\newcommand{\ecnj}{\bbox\end{conjecture}}
\newcommand{\eprb}{\bbox\end{problem}}
\newcommand{\eque}{\bbox\end{question}}
\newcommand{\ecor}{\end{corollary}}
\newcommand{\efac}{\end{fact}}
\newcommand{\eprf}{\bbox}
\newcommand{\beqn}{\begin{equation}}
\newcommand{\eeqn}{\end{equation}}
\newcommand{\wbox}{\mbox{$\sqcap$\llap{$\sqcup$}}}
\newcommand{\bbox}{\vrule height7pt width4pt depth1pt}
\newcommand{\qed}{\bbox}

\newcommand{\rarrow}{\rightarrow}
\newcommand{\larrow}{\leftarrow}
\newcommand{\grad}{\bigtriangledown}
\overfullrule=0pt
\def\setof#1{\lbrace #1 \rbrace}

\section{The permanent}
Recall that a $n \times n$ matrix $A$ is called doubly stochastic if it is nonnegative entry-wise
and its every column and row sum to one. The set of $n \times n$ doubly stochastic
matrices is denoted by $\Omega_{n}$. The set of $n \times n$ of row stochastic(i.e. when every row sum to one)
is denoted by $RS_{n}$, the set of column stochastic(i.e. when every column sum to one)
is denoted by $CS_{n}$.

Let $\Lambda(k,n)$ denote the set of
$n \times n$ matrices with nonnegative integer entries and row and column sums all
equal to $k$ . We define the following subset of rational doubly stochastic matrices:
$\Omega_{k,n} = \{ k^{-1} A : A \in \Lambda(k,n) \}$.\\
  
Recall that the permanent of a square matrix A is defined by
$$ 
per(A) = \sum_{\sigma \in S_{n}} \prod^{n}_{i=1} A(i,
\sigma(i)). 
$$
The following inequality was conjectured by B.l. van der Waerden in 1926 and proved independently in 1981 by D.L. Falikman \cite{fal} and G.P. Egorychev \cite{ego}:
\beqn \label{celeb}
\min_{A \in \Omega_{n}} per(A) = \frac{n!}{n^n} =: vdw(n).
\eeqn

\subsection {Schrijver-Valiant Conjecture and (main) Schrijver's permanental inequality}
Define 
$$
\lambda(k,n) = \min \{per(A) : A \in \Omega_{k,n}\} = k^{-n} \min \{per(A) : A \in \Lambda(k,n)\} ; 
$$
$\theta(k) = \lim_{n \rightarrow \infty}(\lambda(k,n))^{\frac{1}{n}}$.\\

It was proved in \cite{schr-val81} (also earlier in \cite{wilf}) that, using our notations, $\theta(k) \leq G(k) =: (\frac{k-1}{k})^{k-1}$ and conjectured
that $\theta(k) = G(k)$. 
Though the case of $k=3$ was proved by M. Voorhoeve in 1979 \cite{vor} , this conjecture was settled
only in 1998 \cite{schr-val} (17 years after the published proof of the Van der Waerden Conjecture). The main result of \cite{schr-val} (as many people, including myself, wrongly
thought) 
is the remarkable {\bf (Schrijver-bound)} :\\
\beqn \label{S-B}
\min \{per(A) : A \in \Omega_{k,n}\} \geq \left(\frac{k-1}{k} \right)^{(k-1)n} 
\eeqn
The bound (\ref{S-B}) is a corollary of another inequality for doubly-stochastic matrices:
\beqn \label{S-M}
per(\widetilde{A}) \geq \prod_{1 \leq i,j \leq n} (1- A(i,j)); A \in \Omega_{n}; \widetilde{A}(i,j) =: A(i,j)(1-A(i,j)), 1 \leq i,j \leq n.
\eeqn

The proof of (\ref{S-M}) in \cite{schr-val} is, in the words of its author, "highly complicated". Surprisingly,
the only known to me application of (\ref{S-M}) is the bound (\ref{S-B}), which applies only to "very" rational
doubly-stochastic matrices. The main goal of this paper is to show the amazing power of (\ref{S-M}), which
has been overlooked for 13 years.

\section{A Generalization of Schrijver's permanental inequality}
We prove in this section the following theorem, stated in \cite{von} in a rather cryptic way.Fortunately, the paper
cites \cite{CW} and M. Chertkov is my colleague in Los Alamos.\\ 
The statement
in the current paper has been communicated to me by Misha Chertkov, to whom I am profoundly grateful.
\dfn \label{trata}
Define for a pair $(P,Q)$ of non-negative matrices the following functional:
\beqn \label{Bethe}
CW(P,Q) =: \sum_{1 \leq i,j \leq n} \log(1- Q(i,j)) (1- Q(i,j)) - \sum_{1 \leq i,j \leq n} Q(i,j) \log\left(\frac{Q(i,j)}{P(i,j)} \right).
\eeqn
{\it (Note that for fixed $P$ the functional $CW(P,Q) = \sum_{1 \leq i,j \leq n} F_{i,j}(Q(i,j))$ and $F_{i,j}(0) =0$.\\
If $Q \in \Omega_{n}$ is doubly-stochastic and $P = Diag(a_1,...,a_n)T Diag(b_1,...,b_n)$ then
\beqn \label{sca}
CW(P,Q) = \sum_{1 \leq i \leq n} \log(a_i b_i) + CW(T,Q).
\eeqn
Therefore, WLOG we can consider only doubly-stochastic matrices $P$.\\ 
The functional $CW(P,Q)$ is concave in $P$ and, rather surprisingly (see the 2011 arxiv version of \cite{von}),  concave in $Q \in \Omega_{n}$.}\\
\edfn
\thm \label{THM}
Let $P$ be non-negative $n \times n$ matrix. If $Per(P) > 0$ then $max_{Q \in \Omega_{n}}CW(P,Q)$ is attained and
\beqn
\log(Per(P)) \geq max_{Q \in \Omega_{n}} CW(P,Q)
\eeqn
(It is assumed that $0^0 = 1$.)\\
An equivalent statement of this theorem is
\beqn \label{glav}
\log(Per(P)) \geq \sum_{1 \leq i,j \leq n} \log(1- Q(i,j)) (1- Q(i,j)) - \sum_{1 \leq i,j \leq n} Q(i,j) \log\left(\frac{Q(i,j)}{P(i,j)} \right): P \geq 0, Q \in \Omega_{n}
\eeqn 
\ethm

\prf
We will prove, to avoid trivial technicalities, just the positive case, i.e when $P(i,j) > 0, 1 \leq i,j \leq n$.\\
We compute first partial derivatives:
\beqn \label{pd}
\frac{\partial}{\partial Q} CW(P, Q) = \{ -2 - \log(1- Q(i,j)) - \log(Q(i,j)) + \log(P(i,j)) : 1 \leq i,j \leq n\}
\eeqn 
In the positive case, i.e. for the fixed positive $P$, the functional $CW(P,Q)$ is bounded and continuous on $ \Omega_{n}$. Therefore
the maximum exists. Let $V \in \Omega_{n}$ be one of argmaximums, i.e. 
$$
CW(P,V) = max_{Q \in \Omega_{n}} CW(P,Q).
$$
Then, after some column/row permutations
$$
V = \left( \begin{array}{cccc}
		  V_{1,1} &0&...&0 \\
		  0&V_{2,2}& 0&...0 \\
		  . &.& .& . \\
		  0 &...& 0 & V_{k,k} \end{array} \right);
$$

$$
P = \left( \begin{array}{cccc}
		  P_{1,1} &.&...&. \\
		  .&P_{2,2}& .&.... \\
		  . &.& .& . \\
		  . &...& . & P_{k,k} \end{array} \right);
$$
The diagonal blocks  $V_{i,i}$ are indecomposable doubly-stochastic $d_i \times d_i$ matrices;\\
$\sum_{1 \leq i \leq k} d_i = n$ and $1 \leq k \leq n$. Clearly,
$$
CW(P,V) = \sum_{1 \leq i \leq k} CW(P_{i,i},V_{i,i}).
$$
As $\log(per(P)) \geq \sum_{1 \leq i \leq k} \log(per(P_{i,i}))$ it is sufficient to prove that 
$$
\log(Per(P_{i,i})) \geq CW(P_{i,i},V_{i,i}); 1 \leq i \leq k.
$$
For blocks of size one, the inequality is trivial: $(1-1)^{1-1} - 1 \log(\frac{1}{a}) = \log(a)$.\\
Consider a (indecomposable) block $V_{i,i}$ of size $d_i \geq 2$ and define its support
$$
Supp( V_{i,i}) = \{(k,l): V_{i,i}(k,l) > 0 \}.
$$
Note that $1 > V_{i,i}(k,l) > 0, (k,l) \in Supp(V_{i,i})$. Consider the following functional
$$
L(W_{i,i}) =: \sum_{(k,l) \in Supp(W_{i,i} } \log(1- W_{i,i}(k,l)) (1-  W_{i,i}(k,l)) - \sum_{(k,l) \in Supp(V_{i,i}} W_{i,i}(k,l) \log\left(\frac{ W_{i,i}(k,l)}{P(i,j)} \right)
$$
defined on compact convex subset of doubly-stochastic matrices which are zero outside of $Supp( P_{i,i})$. We conclude that the functional $L(\dot)$ is differentiable
at $ V_{i,i}$. Note that $L(V_{i,i}) = CW(P_{i,i},V_{i,i}))$. \\
We now can express the local extremality condition not on full $\Omega_{d_i}$ but rather on its
compact convex subset of doubly-stochastic matrices which are zero outside of $Supp( P_{i,i})$.\\
Using (\ref{pd}) and doing standard Lagrange multipliers respect to variables $V_{i,i}(k,l), (k,l) \in Supp( V_{i,i})$, we get that there exists real numbers $(\alpha_k; \beta_l)$
such that
$$
-2 -\log(1-V_{i,i}(k,l)) - \log(V_{i,i}(k,l)) + \log(P_{i,i}(k,l)) = \alpha_k + \beta_l: (k,l) \in Supp( V_{i,i}).
$$
Which gives for some positive numbers $a_k, b_l$ the following scaling: 
\beqn
P_{i,i}(k,l) = a_k b_l V_{i,i}(k,l) (1 - V_{i,i}(k,l)); (k,l) \in Supp( V_{i,i}).
\eeqn
It follows from the definition of the support that
\begin{enumerate}
\item
\beqn \label{raz}
P_{i,i} \geq Diag(a_k)\widetilde{V_{i,i}} Diag(b_l); \widetilde{V_{i,i}}(k,l) = V_{i,i}(k,l) (1 - V_{i,i}(k,l)).
\eeqn
\item
Using the scalability (\ref{sca}) property, we get that
\beqn \label{dva}
CW(P_{i,i},V_{i,i}) = \sum \log(a_k) + \sum \log(b_l) + \sum_{(k,l) \in Supp( V_{i,i})} \log(1 - V_{i,i}(k,l)).
\eeqn
\end{enumerate}
Finally it follows from (\ref{dva}) and Schriver's permanental inequality (\ref{S-M}) that
$$
\log(per(Diag(a_k)\widetilde{V_{i,i}} Diag(b_l)) \geq CW(P_{i,i},V_{i,i});
$$
and that
$$
\log(per(P_{i,i})) \geq \log(per(Diag(a_k)\widetilde{V_{i,i}} Diag(b_l)) \geq CW(P_{i,i},V_{i,i}).
$$
\eprf
\rem
Note that the proof does not use concavity of $CW(P,V)$ in $V \in \Omega_n$.
\erem
\section{Corollaries}
\begin{enumerate}
\item
Schrijver's permanental inequality (\ref{S-M}) is a particular case of (\ref{glav}). Indeed
$$
CW(\widetilde{V}, V) = \sum_{1 \leq i,j \leq n} \log(1 - V(i,j)): V \in \Omega_{n}.
$$
\item
Let $P \in \Omega_{n}$ be doubly-stochastic $n \times n$ matrix. Then
$$
log(per(P)) \geq CW(P,P) =\sum_{1 \leq i,j \leq n} \log(1 - P(i,j)) (1 - P(i,j)).
$$
{\bf We get the following important inequality, perhaps the main observation in this paper}:
\beqn \label{pglav}
\frac{per(P)}{F(P)} \geq 1; F(P) =: \prod_{1 \leq i,j \leq n} \left(1- P(i,j)\right)^{1- P(i,j)}; P \in \Omega_{n}
\eeqn
The lower bound (\ref{pglav}) suggests the importance of the following quantity:
$$
UB(n) =: max_{P \in \Omega_{n}} \frac{per(P)}{F(P)} .
$$
It is easy to show that the limit
$$
UB =: \lim_{n \rightarrow \infty} (UB(n))^{\frac{1}{n}}
$$
exists and $1 \leq UB \leq e$. There is obvious deterministic poly-time algorithm to approximate the permanent of nonnegative matrices
within relative factor $UB(n)$. The current best rate is $e^n$. Therefore proving that $UB < e$ is of
major algorithmic importance.
\rem
All previous lower bounds on the permanent of {\bf doubly-stochastic} matrices $P \in \Omega_n$ depend only on the dimension $n$ and the support of $P$. I.e. the previous
bounds are structural. The beauty (and potential power) of our lower bound (\ref{pglav}) is in its explicit dependence on the entries of $P$. We use (\ref{pglav})
in Section 5 to settle important conjecture on the monomer-dimer entropy.
\erem 

\xam
{\bf I.} Let $P = aJ_n + bI_n, a = \frac{1}{2(n-1)}, b = \frac{n-2}{2(n-1)}$,
i.e. the diagonal $P(i,i) = \frac{1}{2}, 1 \leq i \leq n$ and the off-diagonal entries are equal to $\frac{1}{2(n-1)}$.\\
It is easy to see that for these $(a,b)$:
$$
2^{-n+1} \leq per(aJ_n + bI_n) = n! a^{n} \sum_{0 \leq i \leq n} \frac{1}{i!} \left(\frac{b}{a} \right)^i \leq n! a^{n} exp \left(\frac{b}{a} \right).
$$
Non-difficult calculations show that for this $P \in  \Omega_{n}$
\beqn \label{up}
\frac{per(P)}{F(P)} \approx \left(\sqrt{\frac{e}{2}} \right)^n
\eeqn
{\bf II.}Let $P \in \Omega_{2}  = \frac{1}{2} J_2$ be $2 \times 2$ ``uniform'' doubly-stochastic matrix.
The direct inspection gives that
$$
CW(P,Q) \equiv -2 \log(2) = F(P),  Q \in \Omega_{n}.
$$
Consider now the direct sum $P_{2n} \in \Omega_{2n} = \frac{1}{2} J_2 \oplus...\oplus\frac{1}{2} J_2$.
Then
\beqn
max_{Q \in \Omega_{2n}} CW(P_{2n},Q) = \log(F(P_{2n})) = -2n \log(2).
\eeqn
Therefore in this case
\beqn
\frac{per(P_{2n})}{F(P_{2n})} = 2^{n}.
\eeqn
Which gives the following lower bound on $UB(k)$ for even $k$:
\beqn \label{exact-up}
UB(k) \geq (\sqrt{2})^k.
\eeqn
As $max_{Q \in \Omega_{2n}} CW(P_{2n},Q) = \log(F(P_{2n}))$, this class of matrices also provides a counter-example to 
the non-trivial part of {\bf Conjecture 15} in \cite{von}.\\
{\bf Is the bound (\ref{exact-up}) sharp?}
\exam
\item
Recall the main function from \cite{ejCO8}:
$$
G(x) = \left(\frac{x-1}{x} \right)^{x-1}, x \geq 1.
$$ 
Note that for $P \in \Omega_{n}$ the column product 
\beqn \label{cpr}
CPR_{j}(P) = : \prod_{1 \leq i \leq n} (1- P(i,j))^{1- P(i,j)} \geq G(n).
\eeqn
Define $C_j$ as the number of non-zero entries in the $j$th column then
\beqn \label{cpr-S} 
CPR_{j}(P) = : \prod_{1 \leq i \leq n} (1- P(i,j))^{1- P(i,j)} \geq G(C_j).
\eeqn
The inequality (\ref{cpr}) gives a slightly weaker version of the celebrated Falikman-Egorychev-van der Waerden lower bound (\ref{celeb}):
$$
per(P) \geq \prod_{1 \leq j \leq n} CPR_{j}(P) \geq \left(\frac{n-1}{n} \right)^{n(n-1)}
$$
The inequality (\ref{cpr-S}) gives a non-regular real-valued version of {\bf (Schrijver-bound)}:
\beqn \label{cpr-SS}
per(P) \geq \prod_{1 \leq j \leq n} CPR_{j}(P) \geq \prod_{1 \leq j \leq n} G(C_j)
\eeqn
In the worst case, the author's bound from \cite{ejCO8} is better:
\beqn \label{mybound}
per(P) \geq \prod_{1 \leq j \leq n} G\left( min(j,C_j)\right)
\eeqn
Perhaps, it is true that
\cnj
$$
per(P)\geq \prod_{1 \leq j \leq n} G\left( min(j,EC_j)\right) ?
$$
where the effective real-valued degree $EC_j = G^{-1}(CPR_{j}(P))$.
\ecnj
\end{enumerate}

\section{Some historical remarks}
The column products $CPR_{j}(P) = : \prod_{1 \leq i \leq n} (1- P(i,j))^{1- P(i,j)} \geq G(C_j)$ have appeared in the permanent context before.
Let $P = [a|b,..,|b] \in \Omega_{n}$ be doubly-stochastic matrix with $2$ distinct columns. Then (Proposition 2.2 in \cite{arx05})
\beqn \label{old}
Per(P) \geq CPR(1) vdw(n-1).
\eeqn
Let us recall a few notations from \cite{ejCO8} and \cite{LS10}:\\
\begin{enumerate}
\item
The linear space of homogeneous polynomials with real (complex) coefficients
of degree $n$ and in $m$ variables is denoted $Hom_{R}(m,n)$ ($Hom_{C}(m,n)$).\\
We denote as $Hom_{+}(m,n)$ ($Hom_{++}(n,m)$) the closed convex cone
of polynomials $p \in Hom_{R}(m,n)$ with nonnegative (positive) coefficients.
\item
For a polynomial $p \in Hom_{+}(n,n)$ we define its {\bf Capacity} as
\beqn \label{cap}
Cap(p) = \inf_{x_i > 0, \prod_{1 \leq i \leq n }x_i = 1} p(x_1,\dots,x_n) = 
\inf_{x_i > 0} \frac{p(x_1,\dots,x_n)}{ \prod_{1 \leq i \leq n }x_i }.
\eeqn
\item The following product polynomial is associated with a $n \times n$ matrix $P$:
\beqn
Prod_{P}(x_1,\dots,x_n) =: \prod_{1 \leq i \leq n} \sum_{1 \leq j \leq n} P(i,j)x_j.
\eeqn 
The permanent $per(P)$ is the mixed derivative of the polynomial $Prod_{P}$:
\beqn
per(P)= \frac{\partial^n}{\partial x_1 \partial x_2\dots \partial x_n} Prod_{P}(0).
\eeqn

\item
$$
q_{(j)} =: \frac{\partial}{\partial x_j}Prod_{P}(x_1,\dots,x_n): x_j = 0.
$$
Note that the polynomials $q_{(j)} \in Hom_{+}(n-1,n-1)$\\
For example, $q_{(n)} = \frac{\partial}{\partial x_n}Prod_{P}(x_1,\dots,x_{n-1},0)$.
\end{enumerate}

The following lower bound, which holds for all $P \in \Omega_n$, was proved in \cite{LS10}:
\beqn \label{UNI}
Cap(q_{(j)}) \geq CPR_{j}(P), 1 \leq j \leq n.
\eeqn 
Combining results from \cite{ejCO8} (i.e. $Per(P) \geq vdw(n-1)Cap(q_{(j)}), 1 \leq j \leq n$) and (\ref{UNI}) gives a different version of (\ref{pglav})
\beqn
per(P) \geq \left(\prod_{1 \leq j  \leq n}CPR_{j}(P)\right)^{\frac{1}{n}} vdw(n-1), P \in \Omega_{n}.
\eeqn
Or better
\beqn
per(P) \geq \left(max_{1 \leq j  \leq n}CPR_{j}(P)\right) vdw(n-1), P \in \Omega_{n}.
\eeqn

Perhaps, it is even true that
\cnj \label{cnj-4}
$$
per(P) \geq \prod_{1 \leq j  \leq n} Cap(q_{(j)}), P \in \Omega_n ?
$$
\ecnj
A general, i.e. not doubly-stochastic and not just ``permanental'', version of Conjecture(\ref{cnj-4}) is the following one:
\cnj \label{cnj-5}
Let $p \in Hom_{+}(n,n)$ be {\bf H-Stable}, i.e. $p(z_1,...,z_n) \ne 0$ if the real parts $RE(z_i) > 0, 1 \leq i \leq n$. In
other words, the homogeneous polynomial $p$ does not have roots with positive real parts. Then the following
inequality holds
\beqn \label{new-in}
\frac{\partial^n}{\partial x_1 \partial x_2\dots \partial x_n} p(0) \geq Cap(p)\prod_{1 \leq  j  \leq n} \frac{Cap(q_{(j)})}{Cap(p)}
\eeqn
\ecnj

\section{Some Partial Results Towards the Main Conjecture(s)}
Let us formalize the main new question in the following Conjecture.
\cnj
Let $P \in \Omega_n$ be doubly-stochastic matrix. Is it true that
\begin{enumerate}
\item {\bf ``Optimizational'' Conjecture}
$$
per(P) \leq (\sqrt{2})^n exp(\max_{Q \in \Omega_n} CW(P,Q)).
$$
It will be explained below that {\bf ``Optimizational'' Conjecture} gives provable deterministic polynomial(but not strongly)
algorithm to approximate $per(P)$ with the factor $(\sqrt{2})^n$
\item {\bf Strong Conjecture}
$$
per(P) \leq (\sqrt{2})^n F(P), F(P) =:  \prod_{1 \leq i,j \leq n} \left(1- P(i,j)\right)^{1- P(i,j)}.
$$
{\bf Strong Conjecture} obviously gives deterministic strongly-polynomial algorithm to approximate $per(P)$ with the factor $(\sqrt{2})^n$
\item {\bf Mild Conjecture}
$$
per(P) \leq (\sqrt{2})^n (F(P))^{c}, 
$$
where $0 < c < 1$ is some universal constant. The case $c = \frac{1}{2}$ seems believable.
As $per(P) \geq F(P)$ thus {\bf Mild Conjecture} gives deterministic strongly-polynomial algorithm to approximate $per(P)$ with the factor\\
$(F(P))^{c-1} \leq \approx e^{n(1-c)} < e^n$.
\end{enumerate}
\ecnj
\subsection{Some Basic Properties of $CW(P,Q)$}
The ``odd entropy'' function $OE(p) = p\log(p) - (1-p)\log(1-p), 0 \leq p \leq 1$ is not convex on $[0,1]$. Yet, when lifted to the Simplex it becomes convex. 
This non-obvious result was proved in recent extended version of \cite{von}. We present below a simpler and more general proof.
\dfn
Call a function $f : [0,1] \rightarrow \R$ {\bf simplex-convex} if the functional
$f_{Sim(1)}(p_1,...,p_n) =:f(p_1)+...+f(p_n)$ is convex on the simplex $Sim_n(1) = \{(p_1,...,p_n): p_i \geq 0, 1 \leq i \leq n; \sum_{1 \leq i \leq n} p_i = 1$.
\edfn
Clearly, if $f$ is convex then it is {\bf simplex-convex} as well. We describe a much wider class of {\bf simplex-convex} functions.

We need two simple facts.
\fac \label{bld}
Let $g : [0,1] \rightarrow \R$ be convex function; $g(0) = 0$. Then this function $g$ is super-additive:
$$
g(t_1+...+t_n) \geq g(t_1)+...+g(t_n): t_i \geq 0, 1 \leq i \leq n; t_1+...+t_n \leq 1.
$$
\efac
\prf
Let $t_i \geq 0, 1 \leq i \leq n; t_1+...+t_n = s \leq 1$. Lift $g$ to the simplex $s Sim_n(1) =: Sim_n(s) = \{ t_i \geq 0, 1 \leq i \leq n; t_1+...+t_n = s\}$:
$$
\bar{g}(t_1,...,t_n) =  g(t_1)+...+g(t_n).
$$
The functional $\bar{g}$ is convex on the simplex $Sim_n(s)$. Therefore, its maximium is attained at the extreme points, i.e at the vectors
$(s,0,...,0),...,(0,0,...,s)$. As $g(0) =0$ we get that
$$
\max_{Sim_n(s)} \bar{g}(t_1,...,t_n) = g(s),
$$
which finishes the proof.
\eprf

\fac \label{conv}
Consider the linear subspace $Sum_0 \subset {\bf R}^{n+1}$, $Sum_0 = \{(x_0,x_1,...,x_n): \sum_{0 \leq i \leq n} x_i = 0$.\\
Let $a_0 > 0; a_i > 0, 1 \leq i \leq n$. The diagonal matrix 
$D = Diag(-a_0, a_1,...,a_n)$ is positive semidefinite on $Sum_0$, i.e $<DX,X> \geq 0, X \in Sum_0$
iff
\beqn
\sum_{1 \leq i \leq n} \frac{1}{a_i} \leq  \frac{1}{a_0}; a_i > 0, 1 \leq i \leq n.
\eeqn
\efac
\prf
The proof directly follows from the following easily checkable equality:
$$
\min_{y_1+....+y_n =1} \sum_{1 \leq i \leq n} a_i y_i^2 = \left(\sum_{1 \leq i \leq n} a_i ^{-1} \right)^{-1}.
$$
\eprf

\thm
Let $f : [0,1] \rightarrow \R$ be continuous and twice differentiable on $(0,1)$ function. Define
$g(t) = f(\frac{1}{2} + t),  t \in [-\frac{1}{2},  \frac{1}{2}]$. Assume that the second derivative $g^{(2)}$ satisfies the
following properties:
\begin{enumerate}
\item
$g^{(2)}(t) > 0, 0  >  t > - \frac{1}{2};g^{(2)}(0) \geq 0 $.
\item
$g^{(2)}(t) \geq -g^{(2)}(-t), 0 < t < \frac{1}{2}$.
\item
$\lim_{ t \rightarrow - \frac{1}{2}} \frac{1}{g^{(2)}(t)} = 0$.
\item
The function $\frac{1}{g^{(2)}(t)}$ is convex on $[-\frac{1}{2},0)$.
\end{enumerate}
Then the function $f$ is {\bf simplex-convex}.
\ethm
\rem
The ``odd entropy'' function $OE(p) = p\log(p) - (1-p)\log(1-p), 0 \leq p \leq 1$ satisfies the above properties:\\
Indeed, $\frac{1}{g^{(2)}(t)} = \frac{1}{OE^{(2)}(\frac{1}{2} + t)} = -\frac{1}{8t} + \frac{t}{2}$.
\erem

\prf
As $f$ is continuous it is sufficient to prove that $f(t_0) +...+f(t_n)$ is convex in the interior of the simplex $Sim_{n+1}(1)$, i.e when $ 0 < t_i < 1$.
Define $d_i =: f^{(2)}(t_i)$. We need to prove that $D =: Diag(d_0, d_1,...,d_n)$ is positive semidefinite on $Sum_0$.\\
If $t_i \leq \frac{1}{2}$ then $d_i \geq 0$
and $D$ is positive semidefinite. Otherwise, there is only one $t_i > \frac{1}{2}$, say 
$$
t_0 = \frac{1}{2} + s_0, \frac{1}{2} \geq s_0 > 0; t_i = \frac{1}{2} - s_i; \frac{1}{2} > s_i > 0, 0 \leq i \leq n.
$$
Note that $d_i > 0, 1 \leq i \leq n$. If $f^{(2)}(t_0) = g^{(2)}(s_0) \geq 0$ we are done. Assume that $-\beta =: f^{(2)}(t_0) = g^{(2)}(s_0) < 0, \beta > 0$. Our goal,
using Fact(\ref{conv}), is to prove that
\beqn \label{des}
(\beta)^{-1} \geq \sum_{1 \leq i \leq n} d_i^{-1}.
\eeqn 
Note that
$$
\frac{1}{2} - s_0 = \sum_{1 \leq i \leq n} (\frac{1}{2} - s_i).
$$
Using the properties(2-4) above and Fact(\ref{bld}), applied to the convex function $ \alpha(t) =: \frac{1}{g^{(2)}(t)}, t \in [0,\frac{1}{2}] $, we get that
$$
(f^{(2)}(\frac{1}{2} - s_0))^{-1} \geq \sum_{1 \leq i \leq n} (d_i)^{-1}.
$$
As $f^{(2)}(\frac{1}{2} + s_0) < 0$ we get from property(2) above that
$$
\beta = -f^{(2)}(\frac{1}{2} + s_0) \leq f^{(2)}(\frac{1}{2} - s_0).
$$
Which gives the desired inequality (\ref{des}).
\eprf
\rem
Recall the definition of the {\it Bregman Distance} associated with a convex functional $f$:
$$
0 \leq D_{f}(X||Y) = f(X) - f(Y) - <\grad F_{Y}, X-Y>.
$$
For instance, the {\bf Kullback-Leibler Divergence} is the {\it Bregman Distance} associated with 
$$
f(p_1,...,p_n) = \sum_{1 \leq i \leq n} p_i \log(p_i).
$$
As we know that the ``odd entropy'' functional 
$$
OE(p_1,...,p_n)) = \sum_{1 \leq i \leq m} p_i \log(p_i) - (1-p_i)\log(1-p_i)  
$$
is convex on the simplex $Sim_{n}(1)$, we can define a new divergence, which we call {\bf Bethe Divergence}:
\beqn \label{bdiv}
BD(X||Y) = \sum_{1 \leq i \leq n} \left(x_i \log(\frac{x_i}{y_i}) - (1-x_i) \log \left(\frac{1-x_i}{1-y_i} \right) \right); X,Y \in Sim_{n}(1).
\eeqn
It would be interesting to investigate statistical (or learning) applications of the {\bf Bethe Divergence}.
\erem
\subsection{Some easy exact computations of $\max_{Q \in \Omega_n} CW(P,Q)$}
The following fact is easy corollary of the {\bf simplex-convexity} of the ``odd entropy'' function $OE(p) = p\log(p) - (1-p)\log(1-p), 0 \leq p \leq 1$.
\fac \label{gran}
\begin{enumerate}
\item 
Let $p_i > 0, 1 \leq i \leq n$ be a positive vector, $n \geq 2$. Define
$$
OD(q,p) = \sum_{1 \leq i \leq n} (1-q_i)\log(1-q_i) - q_i \log(\frac{q_i}{p_i}), q \in  Sim_n(1).
$$
Then
$$
\max_{(q_1,...,q_n) \in Sim_n(1)} OD(q,p) = \log(p_{j})
$$
iff $p_j \geq \sum_{i \neq j} p_i$. We call such index $j$ dominant.\\
Note that if $n \geq 3$ then there exists at most one dominant index.\\
If there is no dominant index then the maximum is attained in the interior of the simplex $Sim_{n}(1)$.
\item
Let $p_i = const > 0,1 \leq i \leq n$ . Then
$$
\max_{(q_1,...,q_n) \in Sim_n(1)} OD(q,p) = OD(\frac{e}{n},p) = (n-1) \log(1 - n^{-1}) + \log(n) + \log(const).
$$  
\end{enumerate}
\efac

\rem
The first item of Fact (\ref{gran}) says that for $n \geq 3$ the extremum of $OD(q,p), p > 0$ is either an extreme point of the simplex(when the unique dominant
index exists) or a point in the interior. This is in stark contrast with $KLD$-minimization, where the extremum has largest possible support.
\erem
We will take advantage of the following corollary.
\cor \label{gugu}
Let $RS_n$ denote the set of $n \times n$ row-stochastic matrices.\\
Let $P$ be $n \times n$ diagonally dominant non-negative matrix. i.e. $P(i,i) \geq  \sum_{j \neq i} P(i,j); 1 \leq i \leq n$. 
Then
\beqn \label{dom}
\max_{Q \in \Omega_n} CW(P,Q) = \max_{Q \in RS_n} CW(P,Q) = \sum_{1 \leq i \leq n} \log(P(i,i)).
\eeqn
\ecor

The following observation follows now from the scalability property (\ref{sca}). 
\cor \label{gugu1}
Assume that there exist two diagonal matrices $D_1, D_2$ such that the matrix
$P = D_1 A D_2$ is diagonally dominant, i.e. $B(i,i) \geq \sum_{j \neq i} B(i,j)$. Then
$$
\max_{Q \in \Omega_n} CW(P,Q) = \sum_{1 \leq i \leq n} \log(P(i,i)).
$$

\ecor

\subsection{Regular Bipartite Graphs}
Let $RB(r,n)$ denote the set of
$n \times n$ boolean matrices with row and column sums all
equal to $r$. Note that if $A \in RB(r,n)$ then $\frac{1}{r} A$ is doubly-stochastic and
$$
F(\frac{1}{r} A) = \left(\frac{r-1}{r} \right)^{n(r-1)} = G(r)^n.
$$
The celebrated Bregman's upper bound \cite{breg} gives that
$$
per(\frac{1}{r} A) \leq \left( \frac{r!^{\frac{1}{r}}}{r} \right)^n =: B_r^{n}.
$$
Therefore
$$
\frac{per(\frac{1}{r} A)}{F(\frac{1}{r} A)} \leq \left(\frac{B_r}{G(r)}\right)^n \leq \left(\frac{B_2}{G(2)}\right)^n = (\sqrt{2})^n.
$$
Therefore, {\bf Strong Conjecture} holds on the sets $RB(r,n)$.\\
Let $CO(\frac{1}{r}RB(r,n))$ be the convex hull. It follows
from linearity of the permanent in individual rows that
$$
per(\frac{1}{r} A) \leq \left( \frac{r!^{\frac{1}{r}}}{r} \right)^n = B_r^{n}, A \in CO(\frac{1}{r}RB(r,n)).
$$
The following observation(most likely known) follows fairly directly from the classical J.Edmonds' result that the intersection
of convex hulls of incidence vectors of bases of two matroids on the same ground set is equal to the convex hull of incidence vectors of common bases.
\pro
The convex hull
$$
CO(\frac{1}{r}RB(r,n)) = \{A \in \Omega_n: A(i,j) \leq \frac{1}{r}; 1 \leq i,j \leq n \}.
$$
\epro
\cor
$$
CO(\frac{1}{r+1}RB(r+1,n)) \subset CO(\frac{1}{r}RB(r,n)), 1 \leq r \leq n-1.
$$
\ecor
We only can state (rather trivial) upper bound
\beqn \label{rty}
\frac{per(P)}{F(P)} \leq \frac{B_r^{n}}{G(n)^n} \leq \left( \frac{r!^{\frac{1}{r}} e}{r} \right)^n: P \in CO(\frac{1}{r}RB(r,n)).
\eeqn
It follows from (\ref{rty}) that {\bf Strong Conjecture} holds on $CO(\frac{1}{r}RB(r,n)), r \geq 6$.

\subsection{Diagonally Dominant Matrices}
\lem
Let $A$ be $n \times n$ non-negative matrix. Then
\beqn \label{w-had}
Per(A) \leq \prod_{1 \leq i \leq n} (A(i,i)^2 + (\sum_{j \ne i} A(i,j)^2))^{\frac{1}{2}}.
\eeqn
\elem
\prf 
Follows from linearity of the permanent in individual rows and the following generalized
Holder's inequality
\beqn
|\prod_{1 \leq i \leq n} a_i + \prod_{1 \leq i \leq n} b_i| \leq \prod_{1 \leq i \leq n} (|a_i|^n + |b_i|^n)^{\frac{1}{n}}.
\eeqn
\eprf

\cor
If $A$ is Diagonally Dominant then the {\bf ``Optimizational'' Conjecture} holds, i.e.
$$
per(A) \leq (\sqrt{2})^n exp(\max_{Q \in \Omega_n} CW(A,Q)).
$$
\ecor

\section{A proof of Friedland's {\it Asymptotic Lower Matching Conjecture}}
\subsection{Two models for random regular bipartite graphs with multiple edges}
We denote as $RI(r,n)$ the set of $n \times n$ non-negative integer matrices with row and column sums all
equal $r$:
$$
RI(r,n) = \{ \{A(i,j); 1 \leq i,j \leq n\}: A(i,j) \in Z_{+}; Ae = A^{T}e = re\}.
$$
\begin{enumerate}
\item {\bf The Pairing Model}: Consider a random, respect to uniform distribution, permutation $\pi \in S_{rn}$ of length $rn$ and
its standard matrix representation, pictured as a block matrix:\\
$$
M_{\pi} =  \left( \begin{array}{cccc}
		  M_{\pi}(1,1) & M_{\pi}(1,2)&...&M_{\pi}(1,r) \\
		  ... & ... &...&...\\
		  
		  M_{\pi}(r,1) & M_{\pi}(r,2)&...&M_{\pi}(r,r) \end{array} \right),
$$

where each block is a (boolean) $n \times n$ matrix. {\bf The Pairing Model} for a random matrix in $RI(r,n)$ corresponds
to a random matrix $BM(r,n) =: \sum_{1 \leq i,j \leq r} M_{\pi}(i,j)$. This model was used in the context of the permanent in
\cite{schr-val81}.
\item {\bf The sum of $r$ independent permutation matrices}: Another model is just the sum of $r$ independent permutation matrices:
$$
HW(r,n) =: \sum_{1 \leq i \leq r} M_{\sigma_i},
$$
where $\sigma_i \in S_n, 1 \leq i \leq r$ are independent uniformly disributed permutations of length $n$. This model was used
by Herbert Wilf \cite{wilf}. As in \cite{schr-val81}, the main goal and result of \cite{wilf} was the asymptotics of the expected value of the permanent:
\beqn \label{models}
\lim_{n \rightarrow \infty} (E(per( HW(r,n)))^{\frac{1}{n}} = \lim_{n \rightarrow \infty} (E(per( BM(r,n)))^{\frac{1}{n}} = r G(r).
\eeqn
\end{enumerate}
It is worth noticing that the proof in \cite{wilf} is much more involved than in \cite{schr-val81}. One of the corollaries of (\ref{models}) is the
following inequality
\beqn
\lim_{n \rightarrow \infty} (\min_{A \in RI(r,n)}(per(A))^{\frac{1}{n}} \leq r G(r),
\eeqn
which was proved much later to be equality.\\
Let $prob_1(r,n)$ be the probability of the event $BM(r,n) \in RB(r,n)$, where $RB(r,n)$ is the set of  $n \times n$ boolean matrices with $r$ ones in each row and column;
$prob_2(r,n)$ be the probability of the event $HW(r,n) \in RB(r,n)$. Brendan McKay conjectured in \cite{brm-84} that for fixed $r$(we present here a simplified expression)
\beqn \label{popa}
prob_1(r,n) =exp \left(-\frac{(r-1)^2}{2} + O(n^{-1}) \right).
\eeqn
This conjecture was proved almost 20 years after in \cite{brm-03}, moreover it holds for $r = o(\sqrt{n})$. The proof in \cite{brm-84} is rather involved and has nothing to do with the permanent.\\
On the other hand, it is easy to see that
\beqn
\frac{1}{(n!)^{r-1}} \prod_{1 \leq i \leq r-1} \min_{A \in RB(n-i,n)} per(A) \leq prob_2(r,n) \leq \frac{1}{(n!)^{r-1}} \prod_{2 \leq i \leq r} \max_{A \in RB(n-i,n)} per(A).
\eeqn
We can use now various lower bounds on $\min_{A \in RB(n-i,n)} per(A)$ and the Bregman's upper bound $((n-i)!)^{\frac{n}{n-i}}$ on $\max_{A \in RB(n-i,n)} per(A)$.\\
Using just the Van Der Waerden-Falikman-Egorychev bound we get that
\beqn \label{pip1}
\prod_{1 \leq i \leq r-1} (\frac{n-i}{n})^{n} \leq prob_2(r,n) \leq \prod_{1 \leq i \leq r-1} \frac{((n-i)!)^{\frac{i}{n-i}}}{(n-i+1)...n}
\eeqn
The best current lower bound (\ref{mybound}) gives
\beqn \label{pip2}
\prod_{1 \leq i \leq r-1} G(n-i)^{i} \frac{(n-i)^i}{(n-i+1)...n} \leq prob_2(r,n) \leq \prod_{1 \leq i \leq r-1} \frac{((n-i)!)^{\frac{i}{n-i}}}{(n-i+1)...n}
\eeqn
For a fixed $r$, as (\ref{pip1}) as well (\ref{pip2}) give the following asymptotic for $prob_2(r,n)$
\beqn \label{as-ind}
prob_2(r,n) \approx exp \left(-\frac{r(r-1)}{2} \right),
\eeqn
which is less than (\ref{popa}).\\
\rem
The expression (\ref{as-ind}) has interesting probabilistic interpretation. Consider ${r \choose 2}$ events $NOV(i,j), 1 \leq i<j \leq r$, responsible
for non-overlapping of permutations $\sigma_i$ and $\sigma_j$. Then $prob_2(r,n) = prob(\cap_{1 \leq i<j \leq r} NOV(i,j))$ and (\ref{as-ind}) states
that events $NOV((i,j)), 1 \leq i<j \leq r$ are asympotically independent.
\erem
Ian Wanless noticed in \cite{wan} that
$$
\min _{A \in RB(r,n)} per(A) \leq (prob_1(r,n))^{-1} E(per( BM(r,n)).
$$
Together with (\ref{popa}) it implies that
\beqn \label{toot}
\lim_{n \rightarrow \infty} (\min_{A \in RB(r,n)}(per(A))^{\frac{1}{n}}) \leq r G(r),
\eeqn
which is the main conclusion of \cite{wan}. We sketched above an alternative, simpler way to get the same result by combining
Herbert Wilf's 1966 paper and Van Der Waerden-Falikman-Egorychev Inequality and their recent refinements. 

\subsection{Monomer-Dimer Problem}.

Let $per_{m}(A)$ denote the sum of permanents of all  $m \times m$ submatrices of $A$:
$$
per_{m}(A) =: \sum_{|S| = |T| = m} per(A_{S,T}).
$$
Define the following two quantities
$$
EMD_{1}(r,n;m) = E(per_{m}(BM(r,n))), EMD_{2}(r,n;m) = E(per_{m}(HW(r,n))).
$$
A rather direct generalization of derivations in \cite{schr-val81} and \cite{wilf} gives the following asymptotics
$$
\lim_{n \rightarrow \infty, \frac{m}{n} \rightarrow t \in [0,1]} \frac{\log(EMD_{1}(r,n;m))}{n} = g_{r}(t) =: t \log(\frac{r}{t}) -2(1-t)\log(1-t) + (r-t) \log(1 - \frac{t}{r}),
$$
and
$$
\lim_{n \rightarrow \infty, \frac{m}{n} \rightarrow t \in [0,1]} \frac{\log(EMD_{2}(r,n;m))}{n} = g_{r}(t).
$$
It follows that 
\beqn \label{odin}
\lim_{n \rightarrow \infty, \frac{m}{n} \rightarrow t \in [0,1]} \frac{\min_{A \in RI(r,n)} \log(per_{m}(A))}{n} \leq g_{r}(t).
\eeqn
The Wanless argument gives the same inequality for the boolean case
\beqn \label{odinn}
\lim_{n \rightarrow \infty, \frac{m}{n} \rightarrow t \in [0,1]} \frac{\min_{A in RB(r,n)} \log(per_{m}(A)}{n} \leq g_{r}(t).
\eeqn
The Friedland's {\it Asymptotic Lower Matching Conjecture} asserts (after \cite{schr-val81}, \cite{schr-val}) that 
\beqn \label{tup}
\lim_{n \rightarrow \infty, \frac{m}{n} \rightarrow t \in [0,1]} \frac{\min_{A \in RB(r,n)} \log(per_{m}(A))}{n} \geq g_{r}(t).
\eeqn
We prove in this paper a slightly stronger result:
\beqn \label{tup1}
\lim_{n \rightarrow \infty, \frac{m}{n} \rightarrow t \in [0,1]} \frac{\min_{A in RI(r,n)} \log(per_{m}(A)}{n} \geq g_{r}(t).
\eeqn
Of course, as we explained above using Wanless argument, the inequalities ($\geq$) in (\ref{tup}, \ref{tup1}) imply equalities.\\
The {\it Lower Matching Conjecture} asserts that
\beqn
per_{m}(A) \geq D(r;m,n) =: {n \choose m}^2 (\frac{r-t}{r})^{n(r-t)} (tr)^t; A \in RB(r,n), t =: \frac{m}{n}.
\eeqn
We prove in this paper the wollowing weeker inequality but for more general class of matrices, i.e for $A \in RI(r,n)$:
\beqn \label{pa1}
per_{m}(A) \geq SF(r,n,m) =: \frac{(\frac{r-t}{r})^{n(r-t)} (1 - n^{-1})^{(1 - n^{-1})2n^2(1-t)}}{(\frac{t}{r})^{nt} n^{-2n(1-t)} ((n(1-t))!)^2}
\eeqn
We note that
\beqn
\frac{D(r;m,n)}{SF(r,n,m)} = \left(\frac{G(n)^{n-m}}{G(m+1)...G(n)} \right)^2 > 1, m < n,
\eeqn
where $G(x) = (\frac{x-1}{x})^{x-1}, x \geq 1$.\\

The following simple Fact will be used below.
\fac \label{nuzh}
\begin{enumerate}
\item
Define the following function $G(x,t) =(\frac{x-t}{x})^{x-t}, x \geq t \geq 0$.
For a fixed $t > 0$ the function $G(x,t)$ is decreasing in $x$.
\item
Let $(a_1,...,a_k$ be positive numbers, $\sum_{1 \leq i \leq k} a_i = 1$. Then
$$
\prod_{1 \leq i \leq k} (1 - t a_i)^{1 - t a_i} \geq G(k,t); 0 \leq t \leq 1.
$$
\end{enumerate}
\efac

\thm
Let $A \in RI(r,n)$.
For a positive integer $m \leq n$ define 
$$
t = \frac{m}{n}, \alpha = \frac{t}{r}.
$$
Then the following lower bound holds:
\beqn \label{pa}
per_{m}(A) \geq SF(r,n,m) =: \frac{(1-\alpha)^{(1 - \alpha)nr} (1 - n^{-1})^{(1 - n^{-1})2n^2(1-t)}}{\alpha^{nt} n^{-2n(1-t)} ((n(1-t))!)^2}
\eeqn
(Notice that $(1-\alpha)^{(1 - \alpha)nr} = G(r,t)^n$.)
\ethm
\prf
\begin{enumerate}
\item Step 1.\\
Consider the following $2n-m \times 2n-m$ matrix
\beqn \label{kaka}
K =  \left( \begin{array}{cc}
		  aA & bJ_{n,n-m} \\
		  (bJ_{n,n-m})^{T}&0 \end{array} \right),
\eeqn		  
where $a = \alpha = \frac{t}{r}$, $b = \frac{1}{n}$, and $J_{n,n-m}$ is $n \times n-m$ matrix of all ones.\\
It is easy to check that this matrix $K$ is doubly-stochastic. Importantly, the following identity holds:
\beqn \label{ba}
per_{m}(A) = \frac{per(K)}{a^m b^{2(n-m)} ((n-m)!)^2}.
\eeqn 
\item Step 2.\\
We apply the inequality (\ref{pglav}) to the doubly-stochastic matrix $K$\\
\beqn \label{vva}
per(K) \geq \left(\prod_{1 \leq i,j \leq n} (1 - \frac{t}{r} A(i,j))^{(1 - \frac{t}{r})A(i,j)} \right) \left(1 - \frac{1}{n} \right)^{(1 - \frac{1}{n})2n^2(1-t)}.
\eeqn
\item Step 3.\\
Let $d_j$ be the number of non-zero entries in the $j$th column of $A$. Notice that $d_j \leq r$ and $\sum_{A(i,j) \neq 0} \frac{A(i,j)}{r} =1$.
It follows from Fact (\ref{nuzh}) that
\beqn \label{vova}
\prod_{1 \leq i,j \leq n} \left(1 - \frac{t}{r} A(i,j) \right)^{(1 - \frac{t}{r})A(i,j)} \geq \prod_{1 \leq j \leq n} G(d_j,t) \geq G(r,t)^n 
\eeqn 
Which gives the following lower bound on the permanent of $K$:
\beqn \label{voova}
per(K)\geq   G(r,t)^n (1 - \frac{1}{n})^{(1 - \frac{1}{n})2n^2(1-t)}
\eeqn

\item Step 3.\\
Finally, we get (\ref{pa}) by combining the (nontrivial, new) inequality (\ref{voova}) with the (trivial, well known) identity (\ref{ba}).
\end{enumerate}

\eprf
\rem
We can express $SF(r,n,m)$ in terms of the function $G, G(x) = \left(\frac{x-1}{x} \right)^{x-1}, x \geq 1$:
\beqn \label{kuk}
SF(r,n,m) = \frac{(1-\alpha)^{(1 - \alpha)nr}}{\alpha^{nt}} \frac{G(n)^{2n(1-t)}}{(G(1)...G(n(1-t)))^{2} (1-t)^{2n(1-t)}}
\eeqn
\erem
The following more general result is proved in the very same way.
\thm \label{ruru}
Let $P \in \Omega_n$. Then
\beqn \label{pupu}
per_{m}(P) \geq \frac{\left(\prod_{1 \leq i \leq n} (1 - \frac{m}{n}P(i,j))^{1 - \frac{m}{n}P(i,j)} \right) G(n)^{2(n-m)}}{(\frac{m}{n})^m n^{-2n(1-t)} ((n(1-t))!)^2}.
\eeqn 
\ethm
Using Fact(\ref{nuzh}) one can get various corollaries of Theorem(\ref{ruru}) expressed in terms of the support of doubly-stochastic matrix $P$. 

\cor
Fix a positive integer $r$ and consider a sequence of pairs $(n,m)$ such that
$$
n \rightarrow \infty, \frac{m}{n} \rightarrow t \in (0,1).
$$
Then
\beqn \label{kkk}
\frac{\log(SF(r,n,m))}{n} \rightarrow g_{r}(t) = t \log(\frac{r}{t}) -2(1-t)\log(1-t) + (r-t) \log(1 - \frac{t}{r})  
\eeqn
Together with inequalities (\ref{odin}, \ref{odinn}) this solves {\it Asymptotic Lower Matching Conjecture}
\beqn \label{shmu}
\lim_{n \rightarrow \infty} \frac{\log(\min_{A \in RB(r,n)} per_{m}(A))}{n} = \lim_{n \rightarrow \infty} \frac{\log(\min_{A \in RI(r,n)} per_{m}(A))}{n} = g_{r}(t).
\eeqn
\ecor
\prf
We only need to prove (\ref{kkk}). The proof follows either from the Stirling approximation of the factorial or from the representation (\ref{kuk}), using
the well known fact that $\lim_{n \rightarrow \infty} G(n) = e^{-1}$.\\
\eprf

\rem
\begin{enumerate}
\item 
The representation $\frac{n!}{n^n} = \prod_{1 \leq i \leq n} G(i)$ provides very simple derivation of the Stirling formula.
\item
The first published statement of {\it Asymptotic Lower Matching Conjecture} appeared in \cite{fkm}.The author
learned about the statement of (\ref{shmu}) from Shmuel Friedland in 2005.\\ 
The main result of \cite{fried} (and of 2006 arxiv version) was the limit equality (\ref{shmu}) for $t = \frac{r}{r+s}, s = 0,1,2,...$. 
The fairly self-contained and simple 
proof in \cite{fried} was based on the ``hyperbolic polynomials approach'' introduced first in \cite{arx05}. The actual result in \cite{fried}
was stated in terms of {\bf sums of mixed derivatives} of general positive hyperbolic polynomials (the same as {\bf H-Stable} in \cite{ejCO8}),
albeit for a restricted range of the parameter $t$. The proof in the present paper is not general at all, it works only for the $m$-permanent,
i.e. for the class of polynomials $Sym_{m}(y_1,...,y_n)$, where $y_i$ are linear forms with non-negative coefficients. But in this case
the full range of densities $t \in [0,1]$ is covered.\\
Whether it can be generalized to general {\bf H-Stable} polynomials remains open.\\
Our proof of {\it Asymptotic Lower Matching Conjecture} illustrates once more how badly had the ``Bethe Restatement'' of Schrijver's inequality (\ref{S-M})
been overlooked.\\
The author did some search on {\it Google Scholar} and found, to his amazement, that the Bethe approximation is one the oldest heuristics for
the monomer-dimer problem, goes back to 1930s. So, the recent Bethe Approximation approach(as a heuristic) to the permanent is, in a way, a rediscovery.
Apparently, the first recent publication in this direction was \cite{colu}.\\
How cool is it that
this classical statistical physics stuff was one of the main keys to rigorously settle the {\it Asymptotic Lower Matching Conjecture}!
Of course, it would have been rather useless without the amazing Schrijver's
inequality (\ref{S-M}). Note that the validity of Conjecture 4.1 also implies {\it Asymptotic Lower Matching Conjecture}.
It would be great
to prove Conjecture 4.1 using {\bf H-Stable} polynomials.\\
\item
The following equality holds for the doubly-stochastic matrices $K$ as in (\ref{kaka}):
$$
F(K) = \max_{Q \in \Omega_n} CW(K,Q)
$$
\end{enumerate}
\erem
\section{ A disproof of a positive correlation conjecture due to [Lu,Mohr,Szekely]}
Let $A$ be $n \times n$ stochastic matrix, i.e. the rows of $A$ are probabilistic distributions on $\{1,...,n\}$;
$(e_1,...,e_n)$ is the standard basis in $R^n$.\\
Let ${\bf V} =: (V_1,...,V_n)$
be a $n$-tuple of independent random vectors:
$$
Prob(V_i = e_k) = A(i,k); 1 \leq i,k \leq n.
$$ 
The distribution of the sum $V_1+...+V_n$ coincides with the vector of the coefficients of the
product polynomial 
$$
Prod_{A}, Prod_{A}(x_1,...,x_n) = \prod_{1 \leq i \leq n} \sum_{1 \leq j \leq n} A(i,j)x_j,
$$
i.e. the probability $Prob(V_1+...+V_n = (\omega_1,...,\omega_n))$ is the coefficient $a_{\omega_1,...,\omega_n}$ of the monomial
$\prod_{1 \leq i \leq n} x_i^{\omega_i}$ in the polynomial $Prod_{A}$. In particular,
\beqn \label{vert}
per(A) = Prob(V_1+...+V_n = e),
\eeqn
where $e = (1,1,...,1)$ is the vector of all ones.\\
Notice that the expected value $E(V_1+...+V_n) = (c_1,...,c_n)$, where $c_j$ is the sum of the jth column of $A$. Thus in the 
doubly-stochastic case
\beqn \label{conc}
per(A) = Prob(V_1+...+V_n = E(V_1+...+V_n)) = Prob(||V_1+...+V_n - E(V_1+...+V_n)||_{2}^2 < 2),
\eeqn
and the lower bounds on the permanent of doubly-stochastic matrices can be viewed as concentration inequalities
for sums of independent random vectors. This interpretation raises a number of natural questions:
\begin{enumerate}
\item What are the lower bounds on  $Prob(||V_1+...+V_n - E(V_1+...+V_n)||_{2}^2 \leq R), R \leq n(n-1)$ in the doubly-stochastic
case? Van Der Waerden-Falikman-Egorychev gives the lower bound $\frac{n!}{n^n} \approx exp(-n)$ for $R < 2$.\\
This question, albeit for distributions associated with {\bf H-Stable} polynomials, was asked by the author in \cite{newhyp}.
\item Is it possible to use this probabilistic interpretation to get new lower bounds, like (\ref{pglav}) in this paper?
\item The coefficients of the products polynomials $Prod_A, A \geq 0$, and of more general {\bf H-Stable} and {\bf Strongly Log-Concave}
polynomials \cite{newt}, satisfy a lot of log-concave like inequalities. Perhaps they can used to prove new
concentration inequalies of the type we listed above?
\item We invite the reader to raise more questions.
\end{enumerate}
\rem
We presented above very simple and effective ``classical'' generator to sample the distribution 
$Dist = \{a_{\omega_1,...,\omega_n}: (\omega_1,...,\omega_n) \in Z_{+}^n,\omega_1+...+\omega_n = n\}$,
where $a_{\omega_1,...,\omega_n}$ are coefficients of the product polynomial $Prod_{A}$ and $A$ is a stochastic matrix.\\
The similar problem for the doubly-stochastic polynomial
$$
Per_{U}(x_1,...,x_n) = per(U Diag(x_1,...,x_n) U^*),
$$
where $U$ is $n \times n$ complex unitary matrix, is of major importance in {\bf Quantum Computing}. The generator in this paper can be viewed as a classical approximation:
essentially, we approximate $Per_{U}$ by the lower bound $Prod_{B}$, where $B(i,j) = |U(i,j)|^2; 1 \leq i,j \leq n$.\\
If $p \in Hom_{+}(n,n)$ is doubly-stochastic and log-concave on $R_{+}^n$ then its coefficients satisfy the (sharp) inequality
$$
p_{\omega_1,...,\omega_n} \leq \prod_{1 \leq i \leq n} \omega_i^{-\omega_i}.
$$
The permanental polynomials $Per_{U}(x_1,...,x_n)$ have much veaker, yet sharp, upper bounds:
$$
q_{\omega_1,...,\omega_n} \leq \prod_{1 \leq i \leq n} \frac{(\omega_i) !}{\omega_i^{\omega_i}}.
$$
Notice that as the permanental polynomial $Per_{U}$ is doubly-stochastic thus
$$
Per_{U}(x_1,...,x_n) \geq \prod_{1 \leq i \leq n} x_i; x_j \geq 0, 1 \leq j \leq n,
$$
in other words $per(Q) \geq det(Q)$ for PSD matrices $Q \succeq 0$. It is, of course, a well known result due to I. Schur. But our proof
is much simpler and shorter than all previous ones.

\erem

Define the following $n$ events: 
$$
NE_i =\{(V_1,...,V_n): V_i \not \in \{V_j, j \neq i \} \}; 1 \leq i \leq n.
$$
Equivalently
\beqn \label{verot}
per(A) = Prob(\cap_{1 \leq i \leq n} NE_i).
\eeqn
The authors of \cite{szek} noticed that $Prob(NE_i) = \sum_{1 \leq j \leq n} A(i,j) \prod_{k \neq i} (1 - A(k,j)$ and conjectured the following
beautiful positive correlation inequality for doubly-stochastic matrices $A \in \Omega_n$:
\beqn \label{sze}
per(A) \geq LMS(A) =: \prod_{1 \leq i \leq n} Prob(EV_i) = \prod_{1 \leq i \leq n} \sum_{1 \leq j \leq n} A(i,j) \prod_{k \neq i} (1 - A(k,j)?
\eeqn
It is easy to see that $LMS(A) \geq F(A), A \in \Omega_n$ and $LMS(A) = F(A)$ in the regular case, i.e. when $A \in r^{-1} RB(r,n); 1 \leq r \leq n$.
Therefore in this regular case the inequality (\ref{sze}) holds and is equivalent to the {\bf (Schrijver-bound)} (\ref{S-B}).\\ 
Actually, in this regular a stronger correlational inequality follows from (\ref{mybound}):
$$
Prob(\cap_{1 \leq i \leq n} NE_i) \geq \frac{G(1)...G(r)}{G(r)^r} \prod_{1 \leq i \leq n} Prob(EV_i).
$$
Apparently the authors of \cite{szek}
did a substantial numerical validation of the
 conjecture on random matrices of modest size.\\
Surprisingly, the Monomer-Dimer Problem provides a probabilistic counter-example.\\
We will present finite families $F_n \subset \Omega_n$ such
that $LMS(A) = Const, A \in F_n$ but the average with some weigths of the permanent over $F_n$ is exponentially smaller than $Const$.
\rem
\begin{enumerate}
\item
One can ask for a {\bf Sidak-like} \cite{szar} correlational inequality:
\beqn \label{sidak}
per(A) =  Prob(\cap_{1 \leq j \leq n} \{|(V_1+...+V_n)_j -1| < 1\}) \geq \prod_{1 \leq j \leq n} Prob(\{|(V_1+...+V_n)_j -1| < 1\}), A \in \Omega_n ?
\eeqn
It is easy to see that 
$$
Prob(\{|(V_1+...+V_n)_j -1| < 1\}) = \sum_{1 \leq i \leq n} A(i,j) \prod_{m \neq i} (1 - A(m,j).
$$
In the notations of Conjecture (\ref{cnj-4}):
$$
Prob(\{|(V_1+...+V_n)_j -1| < 1\}) = q_{(j)}(1,...,1).
$$
So, the conjectured correlation inequality (\ref{sidak}) can be rewritten as
$$
per(A) \geq SD(A) =: \left(\prod_{1 \leq j \leq n}\sum_{1 \leq i \leq n} A(i,j) \prod_{m \neq i} (1 - A(k,j) \right) = \prod_{1 \leq j \leq n} q_{(j)}(1,...,1), A \in \Omega_n ?
$$
Notice that Conjecture (\ref{cnj-4}) claims a smaller lower bound: instead of $q_{(j)}(1,...,1)$, it uses 
$$
Cap(q_{(j)}) =: \inf_{x_i > 0, \prod_{1 \leq i \leq n }x_i = 1} q_{(j)}(x_1,\dots,x_{j-1},x_{j+1},\dots, x_n).
$$
Similarly to [Lu,Mohr,Szekely] conjecture, (\ref{sidak}) holds in the regular case but fails in general:\\
$SD(A) \geq F(A), A \in \Omega_n$ and $SD(A) = LMS(A) = F(A)$ in the regular case. 

\item 
The Sidak Lemma for the gaussian vectors \cite{szar} plays crucial role in the recent Barvinok's bound \cite{barvi} on the number of perfect matchings in general
regular graphs without small cuts.
\item Is there a direct,i.e probabilistic, way to prove correlational inequalities (\ref{sze}, \ref{sidak}) in the regular case? What makes regular bipartite graphs
so ``correlationally'' special? One possible answer is the following observation.
\pro
Let {\bf H-Stable} polynomial $p \in Hom_{+}(n,n)$ be {\bf $r$-regular}, i.e.
\beqn \label{ddss}
p(e + e_it) = \left(\frac{r+t -1}{r} \right)^r; 1 \leq i \leq n. 
\eeqn
(Notice that {\bf $r$-regular} polynomials are doubly-stochastic, therefore $Cap(p) =1$.)\\
Recall the definition of polynomials $q_{(j)} \in Hom_{+}(n-1,n-1)$:
$$
q_{(j)}(x_k: k \neq j) =: \frac{\partial}{\partial x_j}p(x_1,\dots,x_n): x_j = 0.
$$
Then the polynomials $\left( G(r) \right)^{-1}q_{(j)}$ are doubly-stochastic, $1 \leq j \leq n$. 
\epro
\prf
It is easy to see that a polynomial $q \in Hom_{+}(n-1,n-1)$ is doubly-stochastic iff $Cap(q) \geq 1$ and $q(e) =1$.\\
It follows from (\ref{ddss}) that $\left( G(r) \right)^{-1}q_{(j)}(e) = 1, 1 \leq j \leq n$.\\
Finally, Theorem 4.10 in \cite{ejCO8} 
 gives the lower bound $Cap(\left( G(r) \right)^{-1}q_{(j)}) \geq 1$.
\eprf
\cor
Conjecture (\ref{cnj-4}) holds for {\bf $r$-regular} {\bf H-Stable} polynomials.
\ecor
\end{enumerate}

\erem

\subsection{The Construction}
{\it Let us sketch first our strategy: we know that for dobly-stochastic matrices} $A \in \Omega_n$
$$
LMS(A) \geq F(A); SD(A) \geq F(A); per(A) \geq F(A).
$$
{\it We are after matrices where $F(A)$ is asymptotically close to the permanent, but both $LMS(A)$ and $SD(A)$ are much greater than $F(A)$.\\
It is conceivable that there is an explicit counter-example even for $n=3$. We first prove in this paper, by the probabilistic method, only the existence of such matrices
for large enough $n$. We actually don't need the nontrivial inequality $per(A) \geq F(A), A \in \Omega_n$,
but it has motivated the use of the monomer-dimer problem in the construction below. Secondly, we present a concrete counter-example for $n=135$.}\\

Consider either of two random models in $RI(r,n)$, say a random matrix $BM(r,n) \in RI(r,n)$. In induces a conditional
distribution on $RB(r,n)$, i.e. a random matrix $CBM(r,n) \in RB(r,n)$ with the distribution
$$
Prob(CBM(r,n) = A \in RB(r,n)) = \frac{prob(CBM(r,n) = A \in RB(r,n))}{prob\{BM(r,n) \in RB(r,n)\}}.
$$
The Wanless argument gives that
$$
\lim_{n \rightarrow \infty, \frac{m}{n} \rightarrow t \in [0,1]} \frac{\log(E(per_m(CBM(r,n))))}{n} \leq g_{r}(t).
$$
Let $K \in \Omega_{2n-m}$ be the following random doubly-stochastic matrix  
\beqn \label{kaka1}
K =  \left( \begin{array}{cc}
		  aCBM(r,n) & bJ_{n,n-m} \\
		  (bJ_{n,n-m})^{T}&0 \end{array} \right),
\eeqn
$a = \frac{t}{r}, t = \frac{m}{n}; b = \frac{1}{n}$. By the direct inspection, we get that
$$
LMS(K) = \left(t(1 - \frac{t}{r})^{r-1}(1-\frac{1}{n})^{n(1-t)} + (1-t) (1-\frac{1}{n})^{n-1} \right)^{n} \left((1-\frac{1}{n})^{n(1-t)-1} (1-\frac{t}{r})^{r} \right)^{n(1-t)},
$$
and
$$
F(K) = \left(1 -\frac{t}{r} \right)^{(r-t)n} \left(1-\frac{1}{n} \right)^{(n-1)2n(1-t)},
$$
Recall that 
\beqn \label{ba1}
per_{m}(CBM(r,n)) = \frac{per(K)}{a^m b^{2(n-m)} ((n-m)!)^2}.
\eeqn
The conjecture (\ref{sze}) would imply, if true,
that 
\beqn \label{ba2}
per_{m}(CBM(r,n)) a^m b^{2(n-m)} ((n-m)!)^2 \geq LMS(K).
\eeqn
Which would give
\beqn \label{ba3}
f(r,n,m) =: E(per_{m}(CBM(r,n)) a^m b^{2(n-m)} ((n-m)!)^2) \geq LMS(K).
\eeqn
But
$$
\lim_{n \rightarrow \infty, \frac{m}{n} \rightarrow t \in [0,1]} \frac{\log(f(r,n,m))}{n} = (r-t)\log(1 -\frac{t}{r}) - 2(1-t) =: M_{r}(t).
$$
And\\ 
$\lim_{n \rightarrow \infty, \frac{m}{n} \rightarrow t \in [0,1]} \frac{\log(LMS(K))}{n} = 
\log \left(t(1 -\frac{t}{r})^{r-1}e^{-(1-t)} +(1-t)e^{-1} \right) -(1-t)^2 +$\\
$+ r(1-t)\log \left(1-\frac{t}{r} \right) =: S_{r}(t)$.\\
The final observation is the following strict inequality
$$
S_{r}(t) > M_{r}(t): 0 < t < 1, r \geq 1;
$$
which follows from the strict concavity of the logarithm and the inequality
$$
(1 -\frac{t}{r})^{r-1}e^{-(1-t)} > e^{-1}; 0 < t \leq 1, r \geq 1.
$$

\subsection{A disproof of {\bf Sidak-like} positive correlation conjecture (\ref{sidak})}
By the direct inspection, we get that
$$
SD(K) = \left(1-\frac{1}{n-1} \right)^{n(1-t)}\left((1 -\frac{t}{r})^{r-1} (1-\frac{1}{n-1})^{n(1-t)} \right)^{n}.
$$
Which gives the limit
$$
\lim_{n \rightarrow \infty, \frac{m}{n} \rightarrow t \in [0,1]} \frac{\log(SD(K))}{n} = (r-1)\log(1 -\frac{t}{r}) - 2(1-t) =: L_{r}(t).
$$
As $\log(1 -\frac{t}{r}) < 0$ and  $r-1 < r-t$ for $0< t < 1$, it follows that
$$
L_{r}(t) >  M_{r}(t), 0 < t < 1.
$$
\subsection{A Concrete Counter-Example}
To produce a concrete counter-example to Conjecture (\ref{sze}) we consider the case $r=1$ and $t = \frac{1}{2}$. In other words, we consider the following
family of doubly-stochastic matrices $K_n \in \Omega_{n + \frac{1}{2}n}$, where $n$ is even: 
\beqn \label{kaka2}
K_n =  \left( \begin{array}{cc}
		 \frac{1}{2} I_{n} & n^{-1}J_{n,\frac{1}{2}n} \\
		  (n^{-1}J_{n,\frac{1}{2}n})^{T}&0 \end{array} \right).
\eeqn
Then, 
$$
per(K_n) =  (\frac{1}{2})^{\frac{1}{2}n} \frac{1}{n^n} ((\frac{1}{2}n)!)^2 per_{\frac{1}{2}n}(I_n) = \frac{n!}{n^n} 2^{-\frac{1}{2}n};
$$
and 
$$
LMS(K_n) = \left(\frac{1}{2}(1 - \frac{1}{n})^{\frac{1}{2}n} + \frac{1}{2}(1 - \frac{1}{n})^{n-1} \right)^{n} \left(\frac{1}{2} (1 - \frac{1}{n})^{\frac{1}{2}n - 1} \right)^{\frac{1}{2}n}.
$$
We already know that for $n$ large enough $LMS(K_n) > per(K_n)$. Surprisingly, the smallest such $n = 90$. Which gives $135 \times 135$ counter-example to Conjecture (\ref{sze}).
\rem
The value $t = \frac{1}{2}$ is not optimal, we consider it just to simplify the calculations. The optimal value should be the argmaximum of $S_{1}(t)-M_{1}(t), 0 \leq t \leq 1$,
which is $t \approx 0.721$.\\
Note that $\lim_{r \rightarrow \infty} M_{r}(t), S_{r}(t),L_{r}(t) = t-2, 0 \leq t \leq 1$.
\erem

\section{Credits and Conclusion}
The Definition (\ref{trata}) apparently has rich and important stat-physics meaning centered around so called {\it Bethe Approximation}.{\it Bethe Approximation} is
also one of the main Heuristics in modern practice of Machine Learning, especially in inference on graphical models 
(it is quite rare for a Heuristic from Machine Learning to have such amazing proof power).\\
Although this stat-physics background was not used
in the current paper, it and its developers deserve a lot of praise: don't forget that many very good mathematicians have completely overlooked seemingly
simple Theorem \ref{THM}. It would be fantastic to have a rigorous and readable proof of Theorem \ref{THM} based on new(age) methods. The author is a bit skeptical at this point:
any such proof would essentially reprove very hard Schrijver's permanental bound. The other avenue is to better understand and possibly to simplify the original Schrijver's proof,
perhaps it has some deep stat-physics meaning.\\
It is possible that one can use higher order approximation(the Bethe Approximation being of order two,
it involves marginals of subsets of cardinality two). Luckily, this order two case is covered by Schrijver's lower bound (\ref{S-M}). The higher order cases
will probably need new lower bounds (involving subpermanents?). It looks like a beginning of a beautiful(and hard) new line of research.\\
Our proof of Friedland's {\bf monomer-dimer entropy} conjecture illustrates the power of Theorem \ref{THM}. Interestingly, {\bf monomer-dimer entropy} is the classical
topic in stat-physics. The author is not a physicist,passionately so, even after 11 years at Los Alamos. Yet, there is a certain justice in the coincidence that some roots of
this paper can be traced back to Hans Bethe...what a great group of creative people worked in New Mexico back then!\\
   
\section{Acknowledgements}
The author acknowledges the support of NSF grant 116143.\\
Without {\it Google Scholar} the author would have not stumbled on \cite{von}.\\
It is my great pleasure to thank Pascal Vontobel for many things, including catching 
quite a few typos in the previous version, but especially for stating Theorem(\ref{THM}).

\end{document}